\documentclass[11pt]{amsart}

\usepackage{latexsym,amssymb,amsmath}
\usepackage[dvips]{graphicx}

\def\cc{{\mathcal C}}

\def\ff{{\mathcal F}}

\def\ss{{\mathcal S}}

\def\ie{i.e.\ }

\def\ffi{\varphi}
\def\eps{\varepsilon}
\def\dst{\displaystyle}

\def\supp{{\mathrm{supp}\,}}

\def\C{{\mathbb{C}}}

\def\N{{\mathbb{N}}}

\def\R{{\mathbb{R}}}
\def\S{{\mathbb{S}}}

\def\d{\,{\mathrm{d}}}
\newcommand{\norm}[1]{{\left\|{#1}\right\|}}
\newcommand{\ent}[1]{{\left[{#1}\right]}}
\newcommand{\abs}[1]{{\left|{#1}\right|}}
\newcommand{\scal}[1]{{\left\langle{#1}\right\rangle}}

\newenvironment{definition}[1][]{\vskip3pt\noindent\sl\textbf{Definition.}\ }{\rm\vskip3pt}
\newenvironment{remark}[1][]{\vskip3pt\noindent\textbf{Remark.}\ }{\rm\vskip3pt}
\newenvironment{example}[1][]{\vskip3pt\noindent\textbf{Example.}\ }{\rm\vskip3pt}

\newtheorem{lemma}{Lemma}[section]
\newtheorem{proposition}[lemma]{Proposition}
\newtheorem{theorem}[lemma]{Theorem}
\newtheorem{corollary}[lemma]{Corollary}

\begin{document}

\title{Strong annihilating pairs for the Fourier-Bessel transform}
\author{Saifallah Ghobber}

\address{S.G.\,: D\'epartement Math\'ematiques\\
Facult\'e des Sciences de Tunis\\
Universit\'e de Tunis El Manar\\
Campus Universitaire\\ 1060 Tunis\\
Tunisie}
\email{Saifallah.Ghobber@math.cnrs.fr}

\author{Philippe Jaming}

\address{P.J. and S.G\,: Universit\'e d'Orl\'eans\\
Facult\'e des Sciences\\
MAPMO - F\'ed\'eration Denis Poisson\\ BP 6759\\ F 45067 Orl\'eans Cedex 2\\
France}

\address{P.J. (current address)\,: Institut de Math\'ematiques de Bordeaux UMR 5251,
Universit\'e Bordeaux 1, cours de la Lib\'eration, F 33405 Talence cedex, France}
\email{Philippe.Jaming@gmail.com}

\begin{abstract}
The aim of this paper is to prove two new uncertainty principles for the Fourier-Bessel transform (or Hankel transform).
The first of these results is an extension of a result of Amrein-Berthier-Benedicks, it states that a non zero function $f$ and its
 Fourier-Bessel transform $\ff_\alpha (f)$ cannot both have support of finite measure. The second result states that
 the supports of $f$ and $\ff_\alpha (f)$ cannot both be $(\eps,\alpha)$-thin, this extending a result
 of Shubin-Vakilian-Wolff. As a side result we prove that the dilation of a $\cc_0$-function are linearly independent.
 We also extend Faris's local uncertainty principle to the Fourier-Bessel transform.
\end{abstract}



\maketitle


\section{Introduction}
The uncertainty principle is an essential restriction in Fourier analysis. Roughly speaking, this principle states that a function and its Fourier transform cannot be simultaneously well concentrated. There are numerous mathematical formulations for this principle as well as extensions to other transforms (e.g. Fourier type transforms on various types of Lie groups, other integral transforms...) and we refer to
the book \cite{HJ} and the surveys \cite{folland}, \cite{BD} for further references. Our aim here is to consider uncertainty principles in which concentration is measured in sense of smallness of the support and when the transform under consideration is the Fourier-Bessel transform (also known as the Hankel transform).
This transform arises as {\it e.g.} a generalization
of the Fourier transform of a radial integrable function on Euclidean $d$-space as well as
from the eigenvalues expansion of a Schr\"{o}dinger operator.

\medskip

Let us now be more precise and describe our results. To do so, we need to introduce some notations.
Throughout this paper, $\alpha$ will be a real number, $\alpha>-1/2$. For $1\leq p <+\infty$,
we denote by $L^p_\alpha(\R^+)$ the Banach space consisting of measurable functions $f$ on $\R^+$ equipped with the norm
$$
\norm{f}_{L^p_\alpha}=\left( \int_0^\infty \abs{f(x)}^p\,\mbox{d}\mu_\alpha(x)  \right)^{1/p},
$$
where $\dst\mbox{d}\mu_\alpha (x)=(2\pi)^{\alpha+1}x^{2\alpha+1}\,\mbox{d}x$.
For $f\in L^1_\alpha(\R^+)$, the Fourier-Bessel (or Hankel) transform is defined by
$$
\ff_\alpha(f) (y)=\int_0^\infty f(x)j_\alpha(2\pi xy)\d\mu_\alpha (x),
$$
where $j_\alpha$ is the 
Bessel function given by
$$
j_\alpha(x)=
\frac{J_\alpha (x)}{x^\alpha}:=
\frac{1}{2^\alpha}
\sum_{n=0}^{\infty}
\frac{(-1)^n}{n!\Gamma(n+\alpha+1)}\left(\frac{x}{2}\right)^{2n}.
$$
Note that $J_\alpha$ is the Bessel function of the first kind and $\Gamma$ is the gamma function. The function
$j_\alpha$ is even and infinitely differentiable (also entire analytic). One may show that the Fourier-Bessel
transform extends to an isometry on $ L^2_\alpha(\R^+)$ \ie
$$
 \|\ff_\alpha(f) \|_{L^2_\alpha}= \|f \|_{L^2_\alpha}.
$$

\medskip

Uncertainty principles for the Fourier-Bessel transform have been considered in various places,
{\it e.g.} \cite{Bo,RV} for a Heisenberg type inequality or \cite{V} for Hardy type uncertainty principles
when concentration is measure in terms of fast decay. We will here concentrate on uncertainty principles where
concentration is measured in terms of smallness of support.
Our first result (Proposition \ref{prop:local}) is a
straightforward extension of Faris's local uncertainty principle to the Fourier-Bessel transform which compares the
$L^2_\alpha$-norm of $\ff_\alpha(f)$ on some set $E$ of finite measure to weighted norms of $f$
(see Proposition \ref{prop:local} for details).

Our main concern here are uncertainty principles of the following type:
{\sl a function and its Fourier-Bessel transform cannot both have small support}.
 In other words we are interested in the following adaptation of a well-known notion from Fourier analysis:
\begin{definition}\ \\
\label{num}
Let $S$, $\Sigma$ be two measurable subsets of $\R^+$. Then
\begin{itemize}
\item $(S,\Sigma)$ is a \emph{weak annihilating pair} if, $\supp f \subset S$  and $\supp \ff_\alpha(f)\subset \Sigma$ implies $f=0$.

\item $(S,\Sigma)$ is called a \emph{strong annihilating pair} if there exists $C=C_\alpha(S,\Sigma)$ such that
\begin{equation}\label{int,str}
\norm{f}_{L^2_\alpha}\leq C \Big(\norm{f}_{L^2_\alpha(S^c)} +\norm{\ff_\alpha(f)}_{L^2_\alpha(\Sigma^c)} \Big),
\end{equation}
\end{itemize}
where $A^c=\R^+\backslash A$. The constant $C_\alpha(S,\Sigma)$  will be called the \emph{$\alpha$-annihilation constant of $(S,\Sigma)$}.
\end{definition}

Of course, every strong annihilating pair is also a weak one.
There are several examples of  the Uncertainty Principle of the form
 \eqref{int,str} for the Euclidean Fourier transform . One of them is the Amrein-Berthier theorem \cite{AB}
 which is a quantitative version of a result due to Benedicks \cite{Be}  showing that a pair of sets of finite measure
 is an annihilating pair. It is interesting to note that, when $f\in L^2(\R^d)$ the optimal estimate of $C$,
 which depends only on measures $|S_d|$ and $|\Sigma_d|$, was obtained by F. Nazarov \cite{Na} ($d=1$), while in higher dimension the question is not fully settled unless either $S$ or $\Sigma$
is convex ({\it see} the second author's paper \cite{JA} for the best result today).
Our first result will be the following adaptation of the Benedicks-Amrein-Berthier uncertainty principle:

\medskip

\noindent{\bf Theorem A.}\\
{\sl  Let $S$, $\Sigma$ be a pair of measurable subsets of $\R^+$ with
$\mu_\alpha(S),\mu_\alpha(\Sigma)<+\infty$. Then the pair $(S, \Sigma)$ is a strong annihilating pair.}

\medskip

We will actually show a slightly stronger result, namely that a pair of sets with finite Lebesgue measure is strongly annihilating. The proof of this theorem is an adaptation of the proof for the Euclidean Fourier transform in \cite{AB}.
In \cite{AB}, the fact that the Fourier transform intertwines translations and modulations plays a key role.
This property is no longer available for the Fourier-Bessel transform but we have been able to replace
translations by dilations. As a side result, we prove that the dilates of a $\cc_0$-function are linearly
independent.

\medskip

Another Uncertainty Principle which is of particular interest to us is the Shubin-Vakilan-Wolff
theorem \cite[Theorem 2.1]{SVW}, where so called $\eps$-thin sets are considered. 
The natural notion of $\eps$-thin sets for the Fourier-Bessel transform is the following:

\begin{definition}\ \\
A set $S \subset \R^+$ will be called \emph{$(\eps,\alpha)$-thin} if, for $0\leq x\leq 1$,
$$
\mu_\alpha \bigl( S\cap [x,x+1]\bigr)\le
\eps \mu_\alpha \bigl([x,x+1]\bigr)
$$
and for $x>1$,
$$
\mu_\alpha \left( S\cap \ent{x,x+\frac{1}{x}}\right)\le
\eps\mu_\alpha \left( \ent{x,x+\frac{1}{x}}\right).
$$
\end{definition}

We adapt the proof of \cite{SVW} to show the following theorem:

\medskip

\noindent{\bf Theorem B.}\\
{\sl If $\eps$ is small enough and $S$ and $\Sigma$ are $(\eps,\alpha)$-thin then
$$
\norm{f}_{L^2_\alpha}\leq C \Big(\norm{f}_{L^2_\alpha(S^c)} +\norm{\ff_\alpha(f)}_{L^2_\alpha(\Sigma^c)} \Big),
$$
where $C$ is a constant that depends only on $\eps$ and $\alpha$.}

\medskip

The structure of the paper is as follows: in the next section we introduce some further notations as
well as some preliminary results. In Section $3$ we prove the local Uncertainty Inequality for the
 Fourier-Bessel transform. Section $4$ is devoted to the proof of our Amrein-Berthier-Benedicks type theorem
and in Section $5$ we conclude with rou Shubin-Vakilan-Wolff type result, Theorem B.

\section{Preliminaries}
\subsection{Generalities}
In this section, we will fix some notations. We will denote by $|x|$ and
$\scal{x,y}$ the usual norm and scalar product on $\R^d$. The unit sphere of $\R^d$ is denoted by $\S^{d-1}$
and we endow it with the (non-normalized) Lebesgue measure $\mbox{d}\sigma$, that is
$r^{d-1}\,\mbox{d}r\,\mbox{d}\sigma(\zeta)$ is the polar decomposition of the Lebesgue measure. 
The Fourier transform is defined for $F\in L^1(\R^d)$ by
$$
\widehat{F}(\xi)=\int_{\R^d}F(x)e^{-2i\pi\scal{x,\xi}}\,\mbox{d}x.
$$
Note that $\|\widehat{F}\|_2=\norm{F}_2$ and the definition of the Fourier
transform is extended from $F\in L^1(\R^d)\cap L^2(\R^d)$ to $L^2(\R^d)$ in the
usual way.  With this normalizations,  if $F(x)=f(|x|)$ is a radial function on
$\R^d$, then $\widehat{F}(\xi)=\ff_{d/2-1}(f)(|\xi|)$.
More generally, if $F(x)=F_k(|x|)H_k(x)$, $H_k$ a spherical harmonic of degree $k$
(so that $F(r\zeta)=r^kF_k(r)H_k(\zeta)$, $r>0,\zeta\in\S^{d-1}$),
then the Funk-Hecke Formula leads to  $\widehat{F}(\xi)=i^k\ff_{d/2+k-1}(F_k)(|\xi|)H_k(\xi)$,
{\it see} \cite[Chapter IV.2]{SW} for details.
%


If $S_d$ is a measurable set in $\R^d$, we will write $|S_d|$ for its Lebesgue measure.

For $\alpha>-1/2$, let us recall the \emph{Poisson representation formula}
$$
j_\alpha(x)=\frac{1}{2^\alpha\Gamma\left(\alpha+\frac{1}{2}\right)\Gamma\left(\frac{1}{2}\right)}
\int_{-1}^1 (1-s^2)^\alpha\cos sx\frac{\mbox{d}s}{\sqrt{1-s^2}}.
$$
Therefore, $j_\alpha$ is bounded with $|j_\alpha(x)|\leq j_\alpha(0)=\dst\frac{1}{2^\alpha\Gamma(\alpha+1)}$.
As a consequence,
\begin{equation}
\label{eq:L1infty}
\norm{\ff_\alpha(f)}_\infty\leq\frac{1}{2^\alpha\Gamma(\alpha+1)}\norm{f}_{L^1_\alpha}.
\end{equation}
Here $\norm{.}_\infty$ is the usual essential supremum norm.

From the well-known asymptotic behavior of the Bessel function, we deduce that there is a constant $\kappa_\alpha$ such that \begin{equation}\label{appr}
|j_\alpha(t)|\leq \kappa_\alpha t^{-\alpha-1/2}.
\end{equation}

Further, $\ff_\alpha$ extends to a unitary operator on $L^2_\alpha$,
$\norm{\ff_\alpha(f)}_{L^2_\alpha}=\norm{f}_{L^2_\alpha}$. Finally,
if $\ff_\alpha(f)\in L^1_\alpha(\R^+)$, the inverse Fourier-Bessel transform,  is defined for almost every $x$ by
$$
f (x) =\int_0^\infty \ff_\alpha(f) (y)j_\alpha(2\pi xy)\,\mbox{d}\mu_\alpha (y).
$$

Finally, if $I$ is an interval, $I=[a,b]\subset\R^+$
then $3I$ is the interval with same center as $I$ and ``triple'' length, $3I=[a-(b-a),b+(b-a)]\cap\R^+$.
A simple computation shows that the measure $\mu_\alpha$ is \emph{doubling}: there exists a contant $C_\alpha$
such that, for every interval $I\subset\R^+$, $\mu_\alpha(3I)\leq C_\alpha\mu_\alpha(I)$.

\subsection{Generalized translations}

Following Levitan \cite{Bm}, for any function $f \in C^2(\R^+)$ we define the generalized Bessel translation operator
$$
T^\alpha_yf(x)=u(x,y), \hspace{0.5cm}x,\;y \in \R^+,
$$
as a solution of the following Cauchy problem:
$$
\left(\frac{\d^2}{\d x^2}+ \frac{2\alpha+1}{x} \frac{\d}{\d x}\right)u(x,y)=\left(\frac{\d^2}{\d y^2}+ \frac{2\alpha+1}{y} \frac{\d}{\d y}\right)u(x,y),
$$
with initial conditions $ u(x,0)=f(x)$ and $\frac{\partial}{\partial x}u(x,0)=0$.
Here $\dst \frac{\d^2}{\d x^2}+ \dfrac{2\alpha+1}{x} \dfrac{\d}{\d x}$ is  the differential Bessel operator.
The solution of this Cauchy problem can be written out in explicit form:
\begin{equation}\label{translate}
T_x^\alpha(f)(y)=\frac{\Gamma(\alpha+1)}{\sqrt{\pi}\Gamma(\alpha+1/2)}\int_0^\pi f(\sqrt{x^2+y^2-2xy\cos \theta})(\sin \theta)^{2\alpha } \d \theta.
\end{equation}
By formula \eqref{translate}, the operator $T_x^\alpha$ can be extended to all functions $f\in L^p_\alpha(\R^+)$.

The operator $T_x^\alpha$ can be also written by the formula
$$
T_x^\alpha(f)(y)=\int_0^\infty f(t) W(x,y,t) \d\mu_\alpha(t),
$$
where $W(x,y,t) \d\mu_\alpha(t)$ is a probability measure and  $W(x,y,t)$ is defined by
$$
W(x,y,t)= \begin{cases}\dst
\frac{ 2^{-3\alpha}  \Gamma(\alpha+1) }{ \sqrt{\pi}\Gamma\left(\alpha+\frac{1}{2}\right) }
\frac{\Delta(x,y,t)^{2\alpha-1} }{(xyt)^{2\alpha}}&\mbox{if } \abs{x-y}<t< x+y\\
0 &\mbox{otherwise}
\end{cases}
$$
where
$$
\Delta(x,y,t)=\bigl((x+y)^2 -t^2\bigr)^{1/2}\bigl(t^2-  (x-y)^2\bigr)^{1/2}
$$
is the area of the triangle with side length $x,y,t$.
Thus for reasonable functions $f,\; g$, we have
\begin{equation}
\int_0^\infty f(y)T_x^\alpha(g)(y)\d\mu_\alpha(y)=\int_0^\infty g(y)T_x^\alpha(f)(y)\d\mu_\alpha(y).
\label{eq:comtrans}
\end{equation}
Further, $W(x,y,t)\,\mbox{d}\mu_\alpha(t)$ is a probability measure, so that
for $p\ge1$,
$|T^\alpha_xf|^p\leq T_x^\alpha|f|^p$ thus
$\norm{T^\alpha_x f}_{L^p_\alpha(\R^+)}\leq\norm{f}_{L^p_\alpha(\R^+)}$.

The Bessel convolution $f\ast_\alpha g$ of two reasonable functions $f, \;g$ is defined by
$$
f\ast_\alpha g (x)= \int_0^\infty f(t)T_x^\alpha(g)(t)\d\mu_\alpha(t).
$$
Then \eqref{eq:comtrans} reads $ f\ast_\alpha g=g\ast_\alpha f$. It is also well known that for $\lambda>0$, $T_x^\alpha j_\alpha(\lambda.)(y)=j_\alpha(\lambda x)j_\alpha(\lambda y)$. Therefore,
$$
\ff_\alpha\bigl(T^\alpha_xf)(y)=j_\alpha(2\pi xy)\ff_\alpha(f)(y)
$$
and
$$
\ff_\alpha(f\ast_\alpha g)(x)= \ff_\alpha(f)(x)\ff_\alpha(g)(x).
$$

Note also that if $f$ is supported in $[0,b]$ then $T_xf$ is supported in $[0,b+x]$.

\subsection{Linear independence of dilates}

In this section we will prove that the dilation of a $\cc_0$-function are linearly independent, this result may
be of independent interest and plays a key role in the proof of Theorem \ref{weak}.
Let us first introduce the dilation operator $\delta_{\lambda}$, $\lambda>0$, defined by:
$$
\delta_{\lambda} f (x) = \frac{1}{\lambda^{\alpha+1}} f\left(\frac{x}{\lambda}\right).
$$
It is interest to notice that $\ff_\alpha \delta_{\lambda}= \delta_{\frac{1}{\lambda}}\ff_\alpha$.

We may now prove the following lemma which is inspired by a similar result in \cite{ER} for translations.
\begin{lemma}\ \\
\label{libre}
Any nonzero continuous function on $[0,+\infty)$ such that $ \lim_{x\rightarrow +\infty} f(x)=0$ has linearly independent dilates.
\end{lemma}
\begin{proof}
Suppose that there are some distinct elements $\lambda_1, \cdots, \lambda_n \in \R^+ \backslash\{0\}$
and scalars $c_1, \cdots, c_n \in \C$ satisfying
\begin{equation}
\label{libre1}
\sum_{k=1}^n c_k f\left(\frac{x}{\lambda_k}\right)=0 .
\end{equation}
Assume towards a contradiction that one of the scalars $c_k$ is non-zero.
Write $x= e^{s}$ and $\frac{1}{\lambda_k}=e^{\mu_k}$ with $\mu_k, s \in \R$.
Then Equation \eqref{libre1} is equivalent to
\begin{equation}
\label{libre2}
\sum_{k=1}^n c_k g(\mu_k+s)=0 ,
\end{equation}
where $g(x) = f(e^x)$ is a continuous bounded function on $\R$ and $\lim_{x\rightarrow +\infty} g(x)=0$. We will denote by $\widehat{g}$ the
distributional Fourier transform of $g$. Note that, as $g$ is bounded,
$\widehat{g}$ is a distribution of order $0$.

The distributional Fourier transform of equation \eqref{libre2} implies
$$
\left(\sum_{k=1}^n c_k e^{2i\pi \mu_k s} \right)\widehat{g}=0.
$$
As $\dst\sum_{k=1}^n c_k e^{2i\pi \mu_k s}$ is an entire function, its zero set is discrete, therefore $\widehat{g}$ has a discrete support.
Assume $s_0 \in \supp \widehat{g}$, and let $\eta >0$ be such that
$]s_0-\eta, s_0+\eta[\cap \supp \widehat{g}= \{s_0\}$.
Let $\varphi \in C^\infty(\R)$ with support in $]s_0-\eta, s_0+\eta[$ and such that $ \varphi=1$ on $]s_0-\eta/2, s_0+\eta/2[$. Then $\widehat{g}\varphi$ is a distribution of order $0$ such that $\supp \widehat{g}\varphi=\{s_0\}$.
It follows that $\widehat{g}\varphi=c\delta_{s_0}$ for some $c\in \C$.
But then $ g\ast \check\varphi= c e^{2i\pi s_0 t}$, where $\check \ffi$
 is the inverse Fourier transform of $\ffi$. As $\check\ffi\in\ss(\R)$,
 one easily checks that $\lim_{t\to+\infty}g(t)=0$ implies that  $ \lim_{x\to+\infty} g\ast\check\ffi(x)=0$,
 thus $c=0$. It follows that $\supp\widehat{g}=\emptyset$ which implies $f\equiv 0$.
\end{proof}
\section{Local Uncertainty Inequalities}
Heisenberg's inequality for the Fourier-Bessel transform has been established in \cite{RV} as follows:
$$
\norm{xf}_{L^2_\alpha}\norm{x \ff_\alpha(f)}_{L^2_\alpha}\geq (\alpha+1)\norm{f}^2_{L^2_\alpha}.
$$
It says that if $f$ is highly localized, then $\ff_\alpha (f)$ cannot be concentrated near a single
point, but it does not preclude $\ff_\alpha (f)$ from being concentrated in a small neighborhood or more widely separated points.
In fact, the latter phenomenon cannot occur either, and it is the object of local uncertainty inequality to make this precise.
The first such inequalities for the Fourier transform were obtained by Faris \cite{Fa}, and they were subsequently sharpened and generalized
by Price \cite{Pr, Prr}. The corresponding result for the Fourier-Bessel transform is given in the following
proposition:

\begin{proposition}\label{prop:local}\
\begin{enumerate}
  \item If $0 <s < \alpha+1$, there is a constant $K=K(s,\alpha)$ such that for every $f\in L_{\alpha}^2(\R^+)$ and
every measurable set $E \subset \R^+$ of finite measure $\mu_\alpha(E)<+\infty$,
\begin{equation}\label{local}
\norm{\ff_\alpha(f)}_{L_\alpha^2(E)} \leq K \Big[\mu_\alpha (E)\Big]^{\frac{s}{2(\alpha+1)}} \norm{x^s f}_{L_\alpha^2},
\end{equation}
  \item If $s>\alpha+1$, there is a constant $K'=K'(s,\alpha)$ such that for every $f\in L_{\alpha}^2(\R^+)$ and
every measurable set $E \subset \R^+$ of finite measure $\mu_\alpha(E)<+\infty$,
\begin{equation}\label{local2}
\norm{\ff_\alpha(f)}_{L_\alpha^2(E)} \leq K' \mu_\alpha (E)^{1/2}
\norm{f}_{L_\alpha^2}^{1-\frac{(\alpha+1)}{s}} \norm{x^s f}_{L_\alpha^2}^{\frac{\alpha+1}{s}}.
\end{equation}
\end{enumerate}
\end{proposition}

\begin{proof}
As for the first part take $r>0$ and let $\chi_r= \chi_{\{x: \; 0\leq x < r \}}$ and $\tilde{\chi_r}=1- \chi_r$. We may then write
$$
\norm{\ff_\alpha(f)}_{L_{\alpha}^2(E)} =\norm{\ff_\alpha(f)\chi_E}_{L_\alpha^2}
\leq \norm{\ff_\alpha(f\chi_r)\chi_E}_{L_\alpha^2}+\norm{\ff_\alpha(f\tilde{\chi_r})\chi_E}_{L_\alpha^2},
$$
hence, it follows from Plancherel's Theorem that
$$
\norm{\ff_\alpha(f)}_{L_{\alpha}^2(E)}
\leq \mu_\alpha (E)^{1/2}\norm{\ff_\alpha(f\chi_r)}_{\infty} + \norm{f\tilde{\chi_r}}_{L_\alpha^2}.
$$
Now
\begin{eqnarray*}
\norm{\ff_\alpha(f \chi_r)}_{\infty}&\leq& \frac{1}{2^\alpha\Gamma(\alpha+1)}\norm{f\chi_r}_{L_\alpha^1}
\leq \frac{\norm{x^{-s}\chi_r}_{L_\alpha^2}}{2^\alpha \Gamma(\alpha+1)} \norm{x^{s}f}_{L_\alpha^2}\\
&=&a_\alpha r^{\alpha+1-s}\norm{x^{s}f}_{L_\alpha^2}
\end{eqnarray*}
with $a_\alpha=\dst\frac{\pi^{(\alpha+1)/2}}{\sqrt{2^\alpha(\alpha+1-s)}\Gamma(\alpha+1)}$.
On the other hand,
$$
\norm{f\tilde{\chi_r}}_{L_\alpha^2}
\leq \norm{x ^{-s}\tilde{\chi_r}}_\infty \norm{x^{s}f}_{L_\alpha^2}
= r^{-s}\norm{x^{s}f}_{L_\alpha^2},
$$
so that
$$
\norm{\ff_\alpha(f)}_{L_{\alpha}^2(E)}
\leq \Big( r^{-s}
+ a_\alpha r^{\alpha+1-s}\mu_\alpha (E)^{1/2}\Big)\norm{x^{s}f}_{L_\alpha^2}.
$$
The desired result is obtained by minimizing the right hand side of that inequality over $r>0$.

As for the second part we write
$$
\norm{\ff_\alpha(f)}^2_{L_\alpha^2(E)}\le \mu_\alpha (E) \norm{\ff_\alpha(f)}^2_{\infty}\le
 \frac{\mu_\alpha (E)}{(2^\alpha \Gamma(\alpha+1))^2}\norm{f}^2_{L_\alpha^1}.
$$
Moreover
$$
\norm{f}^2_{L_\alpha^1}=\left(\int_0^\infty (1+x^{2s})^{1/2}|f(x)|(1+x^{2s})^{-1/2} \d\mu_\alpha(x)\right)^2,
$$
by the Cauchy-Schwartz inequality, we have
\begin{eqnarray*}
\norm{f}^2_{L_\alpha^1}&\le&\left( \int_0^\infty \frac{\d\mu_\alpha(x)}{1+x^{2s}} \right)
 \left( \int _0^\infty (1+x^{2s})|f(x)|^2\d\mu_\alpha(x)\right) \\
 &=&\left( \int_0^\infty \frac{\d\mu_\alpha(x)}{1+x^{2s}}\right)
 \Big[\norm{f}^2_{L_\alpha^2}+\norm{x^sf}^2_{L_\alpha^2} \Big].
\end{eqnarray*}
 Replacing $f(x)$ by $f(rx)$, $r>0$, in the last inequality gives
$$
\norm{f}^2_{L_\alpha^1}\le \left( \int_0^\infty \frac{\d\mu_\alpha(x)}{1+x^{2s}}\right)
 \left[r^{2(\alpha+1)}\norm{f}^2_{L_\alpha^2}+r^{2(\alpha+1-s)}\norm{x^sf}^2_{L_\alpha^2} \right],
$$
the desired result is obtained by minimizing the right hand side of that inequality over $r>0$.

An easy computation shows that this proof gives
$$
\dst K(s,\alpha)=\dst\frac{\alpha+1}{\alpha+1-s}\left[  \frac{a_\alpha(\alpha+1-s)}{s}\right]^{\frac{\alpha+1}{s}}
$$
and $ \dst K'(s,\alpha)=\frac{1}{2^\alpha \Gamma(\alpha+1)}
 \left[\frac{s}{\alpha+1} \Big(\frac{s}{\alpha+1}-1 \Big)^{\frac{\alpha+1-s}{s}}
\times \int_0^\infty \frac{\d\mu_\alpha(x)}{1+x^{2s}}\right]^{1/2}.$
\end{proof}

\section{Pairs of sets of finite measure are strongly annihilating}
In this section we will show that, if $S$ and $\Sigma$ have finite measure, then the pair $(S,\Sigma)$ is strongly annihilating. Before proving the general case, let us first notice that if $\alpha$ is a positive half-integer, this can be obtained  by transferring the result for the Euclidean Fourier transform established in \cite{JA} (\cite{Na} for $d=1$). Indeed there exists $c_d$ such that,
for $S_d,\Sigma_d\subset\R^d$ of finite Lebesgue measure, and $F\in L^2(\R^d)$,
\begin{equation}\label{jam2}
\norm{F}_{L^2(\R^d)}\leq c_d e^{c_d|S_d||\Sigma_d|}
\Big(\norm{F}_{L^2( S_d^c)}
+\|\widehat{F}\|_{L^2(\Sigma_d^c)} \Big).
\end{equation}
If we define $S$ and $\Sigma$ as
$$
S_d=\{x\in \R^d:\; \abs{x} \in S\}\quad\mathrm{and}\quad\Sigma_d=\{x\in \R^d:\; \abs{x} \in \Sigma\},
$$
then for every function $f \in L^2_{d/2-1}(\R^+)$, there exists $c$ such that
\begin{equation}\label{jam1}
\norm{f}_{L_{d/2-1}^2}\leq c e^{c\mu_{d/2-1}(S)\mu_{d/2-1}(\Sigma) }
 \Big(\norm{f}_{L^2_{d/2-1}(S^c)} +\norm{\ff_{d/2-1}(f)}_{L^2_{d/2-1}(\Sigma^c)} \Big).
\end{equation}

\begin{remark}
It is conjectured that the constant $c_d e^{c_d|S_d||\Sigma_d|}$ in \eqref{jam2} may be replaced
by  $c_d e^{c_d(|S_d||\Sigma_d|)^{1/d}}$ even when $S_d,\Sigma_d$ are not radial sets.
\end{remark}

We will now consider the general case where $ \alpha>-1/2$. We will still show that if  $S$ and $\Sigma$ have finite measure then the pair $(S,\Sigma)$ is strongly annihilating. Unfortunately a precise estimate like \eqref{jam1} still eludes us unless $ \mu_\alpha(S) \mu_\alpha(\Sigma)$ is small enough (see Lemma \ref{small}).
 In order to prove that the pair $(S,\Sigma)$ is strongly annihilating, we will use an abstract result for
 \cite[I.1.1.A, page 88]{HJ}, for which we need the following notations.

 We consider a pair of orthogonal projections on $L^2_\alpha(\R^+)$ defined by
\begin{equation*}
E_S f= \chi_S f,\hspace{1cm} \ff_\alpha (F_\Sigma f)= \chi_\Sigma \ff_\alpha (f),
\end{equation*}
where  $S$ and $\Sigma$ are measurable subsets of $\R^+$.

\begin{lemma}\ \\
\label{equi}
Let $S$ and $\Sigma$ be a measurable subsets of $\R^+$. Then the following assertions are equivalent:
\begin{enumerate}
  \item $\norm{F_\Sigma E_S}<1$;
  \item There exists a constant $D(S,\Sigma)$ such that for all $f \in L^2_\alpha(R^+)$ supported in $S$
$$
 \norm{f}_{L^2_\alpha}\leq D(S,\Sigma) \norm{F_{\Sigma^c}f}_{L^2_\alpha};
 $$
\item $(S,\Sigma)$ is a strongly annihilating pair \ie: there exists a constant $C(S,\Sigma)$ such that for all $f\in L^2_\alpha(R^+)$
    $$
 \norm{f}_{L^2_\alpha}\leq C(S,\Sigma) \left(\norm{E_{S^c}f}_{L^2_\alpha}+\norm{F_{\Sigma^c}f}_{L^2_\alpha}\right).
 $$
\end{enumerate}
Moreover one may take $D(S,\Sigma)=\Big(1-\norm{F_\Sigma E_S }\Big)^{-1}$ and $C(S,\Sigma)=1+D(S,\Sigma)$.
\end{lemma}

\begin{proof}
For sake of completeness let us recall the proof of $(1)\Rightarrow(2)\Rightarrow(3)$, which is the only fact needed in this paper.

Suppose $f$ is supported in $S$. Then
$$
\norm{F_\Sigma f}=\norm{F_\Sigma E_Sf}\le \norm{ F_\Sigma E_S}  \norm{f}_{L^2_\alpha}.
$$
It follows that
$$
\norm{F_{\Sigma^c}f}_{L^2_\alpha}
\ge \norm{f}_{L_\alpha^2}-\norm{F_\Sigma   f}_{L^2_\alpha}
\ge \Big(1-\norm{F_\Sigma E_S }\Big)\norm{f}_{L_\alpha^2}.
$$
Hence, if $\norm{F_\Sigma E_S}<1$, then
$$
 \norm{f}_{L^2_\alpha}\leq \Big(1-\norm{F_\Sigma E_S }\Big)^{-1} \norm{F_{\Sigma^c}f}_{L^2_\alpha}.
 $$
Let us now show the second implication . Let $f\in L^2_\alpha(\R^+)$, then
 \begin{eqnarray*}
 \norm{f}_{L_\alpha^2} &\le& \norm{E_Sf}_{L_\alpha^2} + \norm{E_{S^c}f}_{L_\alpha^2}\\
 &\le& D(S,\Sigma)\norm{F_{\Sigma^c}E_Sf}_{L^2_\alpha} + \norm{E_{S^c}f}_{L_\alpha^2}\\
&=& D(S,\Sigma)\norm{F_{\Sigma^c}(f-E_{S^c}f)}_{L^2_\alpha} + \norm{E_{S^c}f}_{L_\alpha^2}\\
&\le&D(S,\Sigma)\norm{F_{\Sigma^c}f}_{L^2_\alpha}+D(S,\Sigma)
\norm{F_{\Sigma^c}E_{S^c}f}_{L^2_\alpha}+ \norm{E_{S^c}f}_{L_\alpha^2}.
\end{eqnarray*}
Since $\norm{F_{\Sigma^c}E_{S^c}f}_{L^2_\alpha}\le \norm{E_{S^c}f}_{L_\alpha^2}$, we obtain
\begin{equation}\label{cst c}
\norm{f}_{L_\alpha^2} \le  \Big( 1+D(S,\Sigma)\Big) \left(\norm{E_{S^c}f}_{L_\alpha^2}
+ \norm{F_{\Sigma^c}f}_{L^2_\alpha}\right),
\end{equation}
as claimed.
\end{proof}

Unfortunately, showing that $\norm{F_\Sigma E_S }<1$ is in general difficult. However, the Hilbert-Schmidt norm  $ \norm{.}_{HS}$ is mush easier to compute. In our case, we have the following lemma:
\begin{lemma}\ \\
\label{small}
Let $S$ and $\Sigma$ be a pair of measurable subsets of $\R^+$ with finite Lebesgue measure. Then
$$
\norm{F_\Sigma E_S}_{HS}\le \kappa_\alpha\sqrt{2\pi|S|\Sigma|}
$$
where $\kappa_\alpha$ is a numerical constant that depends only on $\alpha$ given by \eqref{appr}.

In particular, if  $|S|\Sigma|< \kappa_\alpha^{-2}$, then
 for any $f\in L^2_\alpha(\R^+)$,
\begin{equation}\label{cst1}
 \norm{f}_{L^2_\alpha}\le \left( 1+\frac{1}{1-\kappa_\alpha\sqrt{2\pi|S|\Sigma|}}\right)
 \Big( \norm{E_Sf}_{L^2_\alpha}+\norm{F_{\Sigma^c}f}_{L^2_\alpha}\Big).
\end{equation}
\end{lemma}


\begin{proof}
The second part of the lemma follows immediately from the fact that
 $\norm{F_\Sigma E_S }\le\norm{F_\Sigma E_S }_{HS}$.

Since $|\Sigma|< +\infty$ it follows from \eqref{appr} that, for every $x>0$,
$j_\alpha(2 \pi x\cdot)\chi_\Sigma \in L^2_\alpha(\R^+)$. A straightforward computation shows
that $F_\Sigma E_S$ is an integral operator with kernel
\begin{equation}\label{eq:noyau}
N(x,y)=\ff_\alpha\Big(\chi_\Sigma j_\alpha(2 \pi x\cdot)\Big)(y)\chi_S (x).
\end{equation}
From Plancherel's theorem, we deduce that
\begin{eqnarray*}
\norm{ F_\Sigma E_S}_{HS}^2&=&\int_0^\infty \abs{\chi_S (x)}^2
\left( \int_0^\infty \abs{\ff_\alpha(\chi_\Sigma j_\alpha(2\pi x\cdot))(y)}^2 \d\mu_\alpha(y)\right)\d\mu_\alpha(x)\\
&=&\int_0^\infty \abs{\chi_S (x)}^2
\left( \int_0^\infty \abs{\chi_\Sigma (y)}^2\abs{j_\alpha(2\pi xy)}^2
\d\mu_\alpha(y)\right)\d\mu_\alpha(x)\\
&=&(2\pi)^{2\alpha+2}\int_0^\infty \int_0^\infty \chi_S (x)\chi_\Sigma (y)
\abs{j_\alpha(2\pi xy)}^2(xy)^{2\alpha+1}\,\mbox{d}x\,\mbox{d}y\\
&\leq& 2\pi\kappa_\alpha^2|S||\Sigma|
\end{eqnarray*}
using \eqref{appr}.
\end{proof}

%

Let us now be more general, set $\alpha>-1/2$ and $S$, $\Sigma$ two measurable subsets of finite measure.

\begin{theorem}\ \\
\label{weak}
Let $S$, $\Sigma$ be a pair of measurable subsets of $\R^+$ with $0<\abs{S}, \abs{\Sigma}<+\infty$. Then
the pair $(S, \Sigma)$ is a strong annihilating pair.
\end{theorem}

\begin{remark}
Let $S$ be a measurable subset of $\R^+$. Using H\"older's inequality one easily shows that,
for every $\eps>0$ there is a constant $C=C(\alpha,\eps)$ depending only on $\alpha$ and $\eps$
such that the Lebesgue measure $\abs{S}$ satisfies
\begin{equation}\label{rem}
\abs{S}\le 1+C_\eps\mu_\alpha(S)^{\frac{1}{2\alpha+2}+\eps}.
\end{equation}
In particular, Theorem A from the introduction follows directly from Theorem \ref{weak}.

Note that the proof below will not give any estimate
on the $\alpha$-annihilation constant of $(S,\Sigma)$.
\end{remark}

\begin{proof}
According to \cite[I.1.3.2.A, page 90]{HJ}, if $F_{\Sigma}E_S$ is compact (in particular if $F_{\Sigma}E_S$ is Hilbert-Schmidt), then if $(S, \Sigma)$ is a weak annihilating pair, it is also a strong annihilating pair. Let us now show that if $0<\abs{S}, \abs{\Sigma}<+\infty$, then $(S, \Sigma)$ is a weak annihilating pair.

In order to do so, let us introduce some further notations.
We will write $E_S \cap F_\Sigma$ for the orthogonal projection onto the intersection of the ranges of $E_S$ and $ F_\Sigma$ and we denote by $\mbox{Im}T$ the range of linear operator $T$.

We will need the following elementary fact on Hilbert-Schmidt operators:
\begin{equation*}
\dim (\mbox{Im}E_S \cap\mbox{Im} F_\Sigma)=\norm{E_S \cap F_\Sigma}_{HS}^2  \leq \norm{ F_\Sigma E_S}_{HS}^2.
\end{equation*}
As $S$ and $\Sigma$ have finite measure then according to Lemma \ref{small} we deduce that
\begin{equation}\label{dimfini}
\dst \dim (\mbox{Im}E_S \cap \mbox{Im}F_\Sigma)\le \norm{F_{\Sigma}E_S}^2_{HS}<+\infty.
\end{equation}

Assume towards a contradiction that there exists $f_0\neq0$ such that
$S_0:= \supp f_0$ and $\Sigma_0:= \supp \ff_\alpha (f_0)$
have both finite measure $0<\abs{S_0}, \abs{\Sigma_0}<+\infty$.

Let $S_1$ be a measurable subsets of $\R^+$ of finite Lebesgue measure such that $S_0\subset S_1$. Since for $\lambda>0$,
$$
\dst\abs{S_1\cup \lambda S_0}=\norm { \chi_{\lambda S_0}-\chi_{S_1} }_{L^2(\R^+)}^2+
 \scal{ \chi_{\lambda S_0},\chi_{S_1}}_{L^2(\R^+)},
$$
the function $\lambda \mapsto \abs{S_1\cup \lambda S_0}$ is continuous on $(0,+\infty)$.
From this, one easily deduces that, there exists an infinite sequence of distinct numbers $\dst(\lambda_j)_{j=0}^{+\infty} \subset(0,\infty)$
with $\lambda_0=1$, such that, if we denote by $S=\dst\bigcup_{j=0}^{+\infty} \lambda_jS_0$ and
$\Sigma=\dst\bigcup_{j=0}^{+\infty}  \frac{1}{\lambda_j}\Sigma_0$,
\begin{equation*}\label{sequence}
 \abs{S}< 2\abs{S_0}, \quad\abs{\Sigma}< 2\abs{\Sigma_0}.
\end{equation*}
We next define $f_i= \delta_{\lambda_i}f_0$, so that $\supp f_i= \lambda_i S_0$.
Since $ \ff_\alpha(f_i)=  \delta_{\frac{1}{\lambda_i}}\ff_\alpha (f_0)$, we have
$\supp \ff_\alpha(f_i)= \frac{1}{\lambda_i} \Sigma_0$.

As $\supp\ff_\alpha(f_0)$ has finite measure, $f_0$ is continuous on $\R^+$ and $f_0(x)\to0$ when $x\to+\infty$.
It follows from Lemma \ref{libre} that $(f_i)_{i=0}^\infty$ are linearly independent vectors belonging to $\mbox{Im}E_{S} \cap\mbox{Im} F_{\Sigma}$,
 which contradicts \eqref{dimfini}.
\end{proof}


\begin{corollary}
Let $S$, $\Sigma$ be a pair of measurable subsets of $\R^+$ with $0<\abs{S}, \abs{\Sigma}<+\infty$
and let $S_d=\{x\in\R^d\,:\ |x|\in S\}$, $\Sigma_d=\{\xi\in\R^d\,:\ |\xi|\in\Sigma\}$. Then
the pair $(S_d, \Sigma_d)$ is a weak annihilating pair for the Fourier transform:
if $F\in L^2(\R^d)$ is such that $\supp F\subset S_d$ and $\supp \widehat{F}\subset \Sigma_d$, then $F=0$.
\end{corollary}

\begin{proof}
We may write, for almost all $r>0$
$$
F(r\zeta)=\sum_{k\geq 0}F_k(r)r^kH_k(\zeta)
$$
where $H_k(\zeta)$ is a spherical harmonic polynomial of degree $k$ and the series converges in the
$L^2(\R^d)$ sense. As
$$
F_k(r)r^kH_k(\zeta)=\int_{\S^{d-1}}F(r\xi)Z_k(\xi,\zeta)\,\mbox{d}\sigma(\xi)
$$
with $Z_k$ the zonal polynomial of degree $k$, $F_k$ is supported in $S$.
Moreover, the Funk-Hecke Formula gives
$$
F(r\zeta)=\sum_{k\geq 0}i^k\ff_{d/2+k-1}[F_k](r)r^kH_k(\zeta)
$$
so that $\ff_{d/2+k-1}[F_k](r)$ is supported in $\Sigma$. As $(S,\Sigma)$ is annihilating for $\ff_{d/2+k-1}$
$F_k=0$ for all $k$, thus $F=0$.
\end{proof}

\begin{remark}
We do not know whether $S_d,\Sigma_d$ is a strong annihilating pair. Indeed, the proof above 
appealed to Fourier-Bessel transforms of various exponents. To prove that $(S_d,\Sigma_d)$
is a strong annihilating pair this way, we would need to prove that $(S,\Sigma)$
is a strong annihilating pair for each $\ff_{d/2+k-1}$, $k=0,1,\ldots$, with annihilation constants
$C_{d/2+k-1}(S,\Sigma)$ independent of $k$.

Moreover, let us denote by $\nu_d(r\zeta)=\,\mbox{d}r\,\mbox{d}\sigma(\zeta)$, $r>0$ and $\sigma\in\S^{d-1}$,
which should be compared to the Lebesgue measure $r^{d-1}\,\mbox{d}r\,\mbox{d}\sigma(\zeta)$.
It is also natural to conjecture
that if $S_d,\Sigma_d\subset\R^d$ are such that $\nu_d(S_d),\nu_d(\Sigma_d)<+\infty$ then
$(S_d,\Sigma_d)$ is a weak annihilating pair for the Fourier transform.
\end{remark}

\section{A result on $\eps$-thin sets}

\subsection{$\eps$-thin sets}
Results in this section are inspired by the ones of Shubin-Vakilian-Wolff
who proved in \cite{SVW} that pairs of $\eps$-thin sets are strongly annihilating for the Euclidean Fourier
transform. To be more precise, let $0<\eps<1$ and let us define $\rho(x)=\min(1,|x|^{-1})$. A measurable set
$S\subset\R^d$ is said to be \emph{$\eps$-thin} if, for every $x\in\R^d$,
$\bigl|S\cap B\bigl(x,\rho(x)\bigr)\bigr|\leq\eps \bigl|B\bigl(x,\rho(x)\bigr)\bigr|$.
Then

\medskip

\noindent{\bf Theorem (Shubin-Vakilian-Wolff \cite[Theorem 2.1]{SVW})}.\\
{\sl There exists $\eps_0$ such that, for every $0<\eps<\eps_0$ there is a constant $C=C(\eps)$ such that,
if $S,\Sigma\subset\R^d$ are $\eps$-thin, then, for every $f\in L^2(\R^d)$,}
$$
\norm{f}_{L^2(\R^d)}\leq C\bigl(\norm{f}_{L^2( S^c)}+\|\widehat{f}\|_{L^2( \Sigma^c)}\bigr).
$$

\medskip

We will now adapt this result to the Fourier-Bessel transform. In order to do so, we first need to define an
appropriate notion of $\eps$-thin sets for the measure $\mu_\alpha$. We want that
the notion which we introduce coincides with the notion of $\eps$-thin radial sets when $\alpha=d/2-1$.

Let us write $\cc_{r_1,r_2}=\{x\in\R^d\,:\ r_1\leq |x|\leq r_2\}$.

Now, take $S=\{r\zeta\,: r\in S_0,\zeta\in\S^{d-1}\}$ be a radial subset of $\R^d$ that is $\eps$-thin
and let us see how the fact that $S$ is $\eps$-thin translates on $S_0$.

First, let $r>2$. Let $\{x_j\}_{j\in J}$ be a maximal subset of $\cc_{r,r+1/r}$
such that $|x_j-x_k|\geq\min\bigl(\rho(x_j),\rho(x_k)\bigr)$. Then the $B\bigl(x_j,\rho(x_j)\bigr)$ cover
$\cc_{r,r+1/r}$. Moreover, it is easy to check that, if $y\in B\bigl(x,\rho(x)\bigr)$ then $ C^{-1}\rho(x)\le\rho(y)\le C\rho(x)$.
It follows that there is a constant $C_d\geq 1$ such that the balls $B\bigl(x_j,C_d^{-1}\rho(x_j)\bigr)$ are disjoint. But then
\begin{eqnarray*}
|S\cap\cc_{r,r+1/r}|&\leq&\sum\bigl|S\cap B\bigl(x_j,\rho(x_j)\bigr)\bigr|\leq \varepsilon
\sum\bigl|B\bigl(x_j,\rho(x_j)\bigr)\bigr|\nonumber\\
&\le&K\eps\sum\bigl|B\bigl(x_j,C_d^{-1}\rho(x_j)\bigr)\bigr|\le K\eps|\cc_{r-1/2r,r+2/r}|.
\end{eqnarray*}
This can be rewritten in terms of $\mu_{d/2-1}$ as
$$
\mu_{d/2-1}(S_0\cap[r,r+1/r])
\le K\eps\mu_{d/2-1}([r-1/2r,r+2/r])\le K\eps\mu_{d/2-1}([r-1/r,r+1/r])
$$
since the measure $\mu_\alpha$ is doubling.

A similar argument leads also to
$$
\mu_{d/2-1}(S_0\cap[r,r+1])\le K\eps\mu_{d/2-1}([r,r+1])
$$
for $r\leq 1$, where $K$ is a constant that depend only of $\alpha$.
This leads us to introduce the definition of $(\eps,\alpha)$-thin sets given in the introduction.
For the convenience of the reader, let us recall it:

\begin{definition}\ \\
Let $\eps\in(0,1)$ and $\alpha>-1/2$.
A set $S \subset \R^+$ is $(\eps,\alpha)$-thin if, for $0\leq x\leq 1$,
$$
\mu_\alpha \bigl( S\cap [x,x+1]\bigr)\le
\eps \mu_\alpha \bigl([x,x+1]\bigr)
$$
and for $x\geq 2$,
$$
\mu_\alpha \left( S\cap \ent{x,x+\frac{1}{x}}\right)\le
\eps\mu_\alpha \left( \ent{x,x+\frac{1}{x}}\right).
$$
\end{definition}

We will need the following simple lemma concerning those sets:

\begin{lemma}\ \\
\label{lem:anncov}
Let $\eps\in(0,1)$ and $\alpha>-1/2$ and let $S \subset \R^+$ be $(\eps,\alpha)$-thin. Then,
there is a constant $C$ depending only on $\alpha$ such that,
if $a\geq 1$ and $b-a\geq\dst\frac{1}{a}$ are such that
$$
\mu_\alpha \left( S\cap \ent{a,b}\right)\le
C\eps\mu_\alpha \left(\ent{a,b}\right)
$$
while for $b>1$,
$$
\mu_\alpha \left( S\cap \ent{0,b}\right)\le
C\eps\mu_\alpha \left(\ent{0,b}\right).
$$
\end{lemma}

\begin{proof}
For $a\geq 1$, we define the sequence $(a_j)_{j\geq 0}$ by $a_0=a$ and $\dst a_{j+1}=a_j+\frac{1}{a_j}$. 
It is easily seen that $(a_j)$ is increasing and $a_j\to+\infty$. Thus there exists $n$ such that $a_n\leq b\leq a_{n+1}$. Note that $b\geq a+1/a=a_1$ thus $n\geq 1$. Further
$a_{n+1}=a_n+1/a_n\leq b+1/a\leq b+b-a$ thus $\mu_\alpha([a,a_{n+1}])\leq C_\alpha\mu_\alpha([a,b])$. 
It follows that
\begin{eqnarray*}
\mu_\alpha(S\cap[a,b])&\leq&\sum_{j=0}^n\mu_\alpha(S\cap[a_j,a_{j+1}])
\leq\eps\sum_{j=0}^n\mu_\alpha([a_j,a_{j+1}])\\
&=& \eps\mu_\alpha([a,a_{n+1}])\leq C_\alpha \eps\mu_\alpha([a,b]).
\end{eqnarray*}

On the other hand, if $b>2$ then $b\geq 1+1/1$ so that
$$
\mu_\alpha(S\cap[0,b])=\mu_\alpha(S\cap[0,1])+\mu_\alpha(S\cap[1,b])
\leq \eps\mu_\alpha([0,1])+C_\alpha\eps\mu_\alpha([1,b])\leq (1+C_\alpha)\eps \mu_\alpha([0,b])
$$
according to the first part of the proof. For $1<b\le2$,
$$
\mu_\alpha(S\cap[0,b])\leq \mu_\alpha(S\cap[0,2])\leq \eps\mu_\alpha([0,2])\leq C_\alpha \eps\mu_\alpha([0,b])
$$
which gives the second part of the lemma.
\end{proof}

\begin{remark}
We will need the following computations. If $r/x\leq x$ then
\begin{eqnarray}
\mu_\alpha\left( \ent{x-\frac{r}{x},x+\frac{r}{x}}\right)
&=&(2\pi)^{\alpha+1}\int_{x-r/x}^{x+r/x}t^{2\alpha+1}\,\mbox{d}t\leq
(2\pi)^{\alpha+1}\frac{2r}{x}(x+r/x)^{2\alpha+1}\nonumber\\
&\le& (2^3\pi)^{\alpha+1}r x^{2\alpha}\label{eq:trives1}.
\end{eqnarray}
On the other hand, for $r/x\geq x/2$ a similar computation shows that
\begin{equation}
\label{eq:trives2}
\mu_\alpha \left( \ent{0,x+\frac{r}{x}}\right)
\le (18\pi)^{\alpha+1}(\frac{r}{x})^{2\alpha+2} .
\end{equation}
\end{remark}

\begin{example}
It should be noted that a measurable subset $(\eps,\alpha)$-thin may not be of finite Lebesgue measure.

\rm Let $\eps \in(0,1)$, $k \in \N $ and $ S= \dst\bigcup_{k\ge 10^6} \ent{k, k+\frac{\eps}{c k}}$
so that $\abs{S}=+\infty$. Moreover if the constant $c$ is large enough then $S$ is $(\eps,\alpha)$-thin. Indeed if $\dst S\cap \ent{x,x+\frac{1}{x}} \neq \emptyset$ then there exists an integer $k$ such that $x\approx k$ and
\begin{eqnarray*}
\mu_\alpha \left( S\cap \ent{x,x+\frac{1}{x}}\right) &=& \mu_\alpha \left(\ent{k,k+\frac{\eps}{c k}}\cap \ent{x,x+\frac{1}{x}}\right) \\
 &\le&  \frac{\eps k^{2\alpha}}{c}
  \leq \eps\mu_\alpha \left(\ent{x,x+\frac{1}{x}}\right)
\end{eqnarray*}
if $c$ is large enough.
\end{example}

\subsection{Pairs of $\eps$-thin sets are strongly annihilating}

We are now in position to prove the following Uncertainty Principle in the spirit of
\cite[Theorem 2.1]{SVW}.

\begin{theorem}\ \label{th:strong} \\
Let $\alpha>-1/2$.
There exists $\eps_0$ such that, for every $0<\eps<\eps_0$, there exists a positive constant
$C$ such that if $S$ and $\Sigma$ are $(\eps,\alpha)$-thin sets in $\R^+$ then for any $f \in L_\alpha^2(\R^+)$
\begin{equation}\label{wolff}
\norm{f}_{L_\alpha^2} \le C \Big(\norm{f}_{L_\alpha^2(S^c)} + \|\ff_\alpha(f)\|_{L_\alpha^2(\Sigma^c)} \Big).
\end{equation}
\end{theorem}

\begin{proof}
In this proof, we construct two bounded integral operators $K$ and $L$ such that
$K+L=I$.Moreover $KE_S$ and $F_\Sigma L$ are bounded operators on $L^2_\alpha(\R^+)$ with
 $$
 \norm{K E_S }\le C_1 \sqrt{\eps}, \quad \norm{ F_\Sigma L} \le C_2 \sqrt{\eps}.
 $$
From such a situation, the Uncertainty Principle can be easily derived.
As
 $$
 \norm{F_{\Sigma}E_S}=  \norm{F_{\Sigma}(L+K)E_S }
 \le  \norm{F_{\Sigma}L}+ \norm{KE_S },
 $$
then
$$
  \norm{F_{\Sigma}E_S}\le (C_1+C_2)\sqrt{\eps}.
$$
Now if $\dst \eps<\eps_0= \frac{1}{(C_1+C_2)^2}$, using Lemma \ref{equi}, we obtain the desired result
$$
  \norm{ f}_{L^2_\alpha} \le  \left( 1+\frac{1}{1-\sqrt{\frac{\eps}{\eps_0}}} \right)\Big( \norm{E _{S^c}f}_{L^2_\alpha}
  + \norm{F_{\Sigma^c}f}_{L^2_\alpha}\Big).
$$

Now we will show how to construct a pair of such operators $K$ and $L$ via a Littlewood-Paley type decomposition. To do so,
 we fix a real-valued Schwartz function $\psi_0: \R^+ \rightarrow \R$ with
$0\leq\psi_0\le 1$, $\supp\psi_0\subset [0,2]$ and $\psi_0=1$ on $[0,1]$
and let $\phi=\ff_\alpha(\psi_0)$. Note that $\phi$ is also in the Schwartz class.

Next, for $j\geq1$ an integer, we define $\psi_j$ by $\psi_j(x)=\psi_0(2^{-j}x)-\psi_0(2^{-j+1}x)$
so that $\psi_j(x)=\psi_1(2^{-j+1}x)$. Note that $\norm{\psi_j}_{L_\alpha^1}=2^{2(\alpha+1)(j-1)} \norm{\psi_1}_{L_\alpha^1}$,
$\norm{\psi_j}_\infty \le 1$, $\supp \psi_j \subset [2^{j-1},2^{j+1}]$ for $j \ge 1$
and $\dst\sum_{j=0}^{\infty} \psi_j = 1$.

Finally, for $j \in \N$ we let $\phi_j (x)= 2^{2(\alpha+1)j} \phi (2^{j}x)$. Thus
$\norm{\phi_j}_{L_\alpha^1}=\norm{\phi}_{L_\alpha^1}$,
$ \ff_\alpha(\phi_j)(\xi)= \ff_\alpha(\phi)(2^{-j}\xi)$,
$\supp \ff_\alpha(\phi_j) \subset [0,2^{j+1}]$ and $\ff_\alpha(\phi_j)=1$ on $[0,2^{j}]$.

Define now the operators $K$ and $L$ on $L_\alpha^2(\R^+)$ in the following way:
\begin{equation}\label{eq:k}
K f= \sum_{j=0}^{+\infty} \psi_j (\phi_j \ast_\alpha f)
\end{equation}
and
\begin{equation}\label{eq:l}
L f= \sum_{j=0}^{+\infty} \psi_j (f- \phi_j \ast_\alpha f).
\end{equation}
Note that the series in \eqref{eq:k} and \eqref{eq:l} converge pointwise since they have at most three nonvanishing terms
at a given point. It is also clear that $ Kf+Lf=f$. Further, $K$ is given by an integral kernel:
$$
Kf (x)= \int_0^{+\infty} A(x,y) f(y) \d \mu_\alpha(y)
$$
where
\begin{equation}\label{eq:a}
   A(x,y)=  \sum_{j=0}^{+\infty} \psi_j(x) T^\alpha_y \phi_j (x).
\end{equation}
We also have
$$
\ff_\alpha (Lf) (x)= \int_0^{+\infty} B(x,y) \ff_\alpha (f)(y) \d \mu_\alpha(y)
$$
where
\begin{equation}\label{eq:b}
   B(x,y)=  \sum_{j=0}^{+\infty}T^\alpha_x  \ff_\alpha(\psi_j)(y) \Big(1-\ff_\alpha( \phi_j) (y)\Big).
\end{equation}
Notice that
\begin{eqnarray}\label{identity}
  B(x,y)&=&  \sum_{j=0}^{+\infty}T^\alpha_x  \ff_\alpha(\psi_j)(y) \Big(1-\ff_\alpha( \phi_j) (y)\Big)
= \sum_{j=0}^{+\infty}T^\alpha_x  \ff_\alpha(\psi_j)(y) \sum_{k>j} \psi_k(y)\nonumber\\
  &=& \sum_{k=1}^{+\infty}\psi_k(y) \sum_{j=0}^{k-1}T^\alpha_x  \ff_\alpha(\psi_j)(y) =  \sum_{k=1}^{+\infty}\psi_k(y) T^\alpha_x  \phi _{k-1}(y).
\end{eqnarray}
This has the same shape as $A(y,x)$.

The remaining of the proof consists in two lemmas. We will first show that $K$ and $L$ are bounded. This will then be used to show that
$$
 \norm{K E_S }\le C_1 \sqrt{\eps}, \quad \norm{ F_\Sigma L} \le C_2 \sqrt{\eps},
 $$
if $S$ and $\Sigma$ are $(\eps,\alpha)$-thin.

To show that $K$ and $L$ are a bounded operators on $L^2_\alpha(\R^+)$. It will suffice to prove the following lemma related to Schur's test:
\begin{lemma} \label{borne} \ \\
The kernel $A$ satisfies the following bounds:
\begin{equation}\label{eq:supx}
   \sup_{x} \int_0^{+\infty} \abs{A(x,y)} \d \mu_\alpha(y) \le C
\end{equation}
and
\begin{equation}\label{eq:supy}
   \sup_{y} \int_0^{+\infty} \abs{A(x,y)} \d \mu_\alpha(x)\le C,
\end{equation}
where $C$ is an absolute constant.

The same bound holds for $B$.
\end{lemma}

\begin{proof}[Proof of Lemma \ref{borne}]
Formula \eqref{eq:supx} follows from the fact that for a fixed $x$ the sum in \eqref{eq:a} contains at most three nonvanishing terms,
$\norm{\psi_j}_\infty \le 1$ and $\norm{\phi_j}_{L^1_\alpha} =\norm{\phi}_{L^1_\alpha}$. Therefore,
$$
\sup_{x} \int_0^{+\infty}  \abs{A(x,y)} \d \mu_\alpha(y) \le 3 \norm{\phi}_{L^1_\alpha}.
$$
Fix $y$ and note that there are at most three values of $j$ such that $\mbox{dist}(y,\supp \psi_j)<1$. Call this set of $j$'s $P$.
We have
$$
 \int_0^{+\infty}  \abs{A(x,y)} \d \mu_\alpha(x) \le 3 \norm{\phi}_{L^1_\alpha} +
  \sum_{j \notin P}  \int_0^{+\infty}  \abs{\psi_j(x)}. \abs{T_y^\alpha \phi_j(x)}\d \mu_\alpha(x).
$$
Since $\phi $ is a Schwartz function we have
$$
\phi_j(t) \le C 2^{2(\alpha+1)j} (1+2^j t)^{-6(\alpha+1)}
$$
and, for $t \ge 1$,
$$
\phi_j(t) \le C 2^{-4(\alpha+1)j}.
$$
Let $x\ge 0$ and $j \notin P$ such that $\psi_j(x)\neq 0$. Since
$$
T_y^\alpha \phi_j(x)= \int_{\abs{x-y}}^{x+y} \phi_j(t) W(y,x,t) \d \mu_\alpha(t)
$$
and $t \ge \abs{x-y} \ge 1$ then
$$
\abs{T_y^\alpha \phi_j(x)}\le C 2^{-4(\alpha+1)j} \int_{\abs{x-y}}^{x+y} W(y,x,t) \d \mu_\alpha(t)=C 2^{-4(\alpha+1)j}.
$$
Hence
$$
\sum_{j \notin P}  \int_0^{+\infty}  \abs{\psi_j(x)}. \abs{T_y^\alpha \phi_j(x)}\d \mu_\alpha(x)\le
 C\sum_{j \ge 0} 2^{-4(\alpha+1)j} \norm{\psi_j}_{L^1_\alpha},
$$
from which we deduce
$$
 \sup_{y} \int_0^{+\infty} \abs{A(x,y)} \d \mu_\alpha(x)\le  3 \norm{\phi}_{L^1_\alpha}
 + C\sum_{j \ge 0} 2^{-2(\alpha+1)j}\norm{\psi_1}_{L^1_\alpha}
$$
which completes the proof for $A$. According to \eqref{identity}, $A$ and $B$ have the same "shape", the proof immediately adapts to $B$.
\end{proof}

Using Schur's test, it follows that $K$ and $L$ are bounded operators on $L^2_\alpha$.

Now we will show that there are constants $C_1$, $C_2>0$  such that
\begin{equation*}\label{eq:keps}
    \norm{K E_S f}_{L^2_\alpha} \le C_1 \sqrt{\eps} \norm{ f}_{L^2_\alpha}
\end{equation*}
and
\begin{equation*}\label{eq:leps}
    \norm{ F_\Sigma L  f}_{L^2_\alpha} \le C_2 \sqrt{\eps} \norm{ f}_{L^2_\alpha}.
\end{equation*}
Using again Schur's test, it will suffice to prove the following lemma:
\begin{lemma} \label{borne2} \ \\
If $S$ and $\Sigma$ are $(\eps,\alpha)$-thin sets, then
\begin{equation}\label{eq:supx2}
   \sup_{x} \int_S \abs{A(x,y)} \d \mu_\alpha(y) \le C \eps
\end{equation}
and
\begin{equation}\label{eq:supy2}
   \sup_{y} \int_\Sigma \abs{B(x,y)} \d \mu_\alpha(x)\le C \eps.
\end{equation}
\end{lemma}

\begin{proof}[Proof of Lemma \ref{borne2}]
By identity \eqref{identity} it will suffice to prove \eqref{eq:supx2}. We want to estimate
$$
\int_S \abs{A(x,y)} \d \mu_\alpha(y) \le \sum_{j\ge 0} \int_S \abs{\psi_j(x)} \abs{T^\alpha_y\phi_j(x)} \d \mu_\alpha(y).
$$
There are at most three values of $j$ such that $\psi_j(x) \neq 0$, so it will suffice to prove
\begin{equation}\label{claim}
    \int_S  \abs{T^\alpha_y\phi_j(x)} \d \mu_\alpha(y)\le C \eps.
\end{equation}
Fix $x$ and let $j$ be such that $\psi_j(x) \neq 0$. Then $2^{j-1}\leq x\leq 2^{j+1}$. We will write $C$ for a constant that depends only on $\alpha$ and that may change from line to line.

Let us explain the method of computation when replacing $\phi$ by $\chi_{[0,1]}$. Then $|\phi_j (t)|= 2^{2(\alpha+1)j}\chi_{[0,2^{-j}]}(t).$
Moreover
\begin{eqnarray*}
\abs{T^\alpha_y\phi_j(x)} &\leq &\int_0^\pi |\phi_j (\sqrt{x^2+y^2-2xy \cos \theta})| (\sin\theta)^{2\alpha} \d\theta\\
&\leq&2^{2(\alpha+1)j} \int_0^\pi \chi_{[0,2^{-j}]}(\sqrt{x^2+y^2-2xy \cos \theta})| (\sin\theta)^{2\alpha}\d\theta.
\end{eqnarray*}
Note that if
$$
x^2+y^2-2xy \cos \theta= (x-y)^2+4xy \sin^2 (\theta/2)\leq 2^{-2j},
$$
then
$$
|x-y|\leq 2^{-j} \ \ \mbox{and}\ \ |\theta| \leq 2^{-2 (j+1)}.
$$
Therefore
$$
\abs{T^\alpha_y\phi_j(x)}\leq 2^{2(\alpha+1)j}\int_0^{2^{-2(j+1)}} \theta^{2\alpha}\d\theta \le C 2^{2(\alpha+1)j}\times2^{-2(2\alpha+1)(j+1)}\le C2^{-2\alpha j}.
$$
It follows that
\begin{eqnarray*}
   \int_{S}\abs{T^\alpha_y\phi_j(x)} \d \mu_\alpha(y)&=&\int_{S\cap [x-\frac{1}{x},x+\frac{1}{x}]}\abs{T^\alpha_y\phi_j(x)} \d \mu_\alpha(y)\\
 &\le& C 2^{-2\alpha j}\mu_\alpha(S\cap [x-\frac{1}{x},x+\frac{1}{x}])
 \le C  \eps 2^{-2\alpha j}x^{2\alpha}
 \le C\eps.
\end{eqnarray*}
As $\phi$ is a Schwartz function, then
$$
\phi (t)\le C_N \sum _{k\ge 0} 2^{-kN} \chi_{[0,2^{k}]}(t),
$$
 where $N$ is a large integer. Then
$$
|T^\alpha_y\phi_j(x)|\leq C_N 2^{2(\alpha+1)j}\sum _{k\ge 0} 2^{-kN}\int_0^\pi \chi_{[0, 2^{k-j}]}(\sqrt{x^2+y^2-2xy \cos \theta})(\sin\theta)^{2\alpha}\d\theta.
$$
Now, if the integral on the right hand side is non zero, then
 $x^2+y^2-2xy \cos \theta:=(x-y)^2+4xy \sin^2 (\theta/2)\leq 2^{2k-2j}$. This implies that $|x-y|\leq 2^{k-j}$ \ie
$ y \in [x-\frac{2^k}{x},x+\frac{2^k}{x}]\cap \R^+.$
Further, for $k<j$ we also have $|\theta| \leq C 2^{k-2j}$ so that
$$
\int_0^\pi \chi_{[0, 2^{k-j}]}(\sqrt{x^2+y^2-2xy \cos \theta})(\sin\theta)^{2\alpha}\d\theta
\leq \int_0^{C2^{k-2j}} \theta^{2\alpha}\d\theta\leq
C 2^{(2\alpha+1)(k-2j)}.
$$
For $k\geq j$, we will use the straightforward inequality
$$
\int_0^\pi \chi_{[0, 2^{k-j}]}(\sqrt{x^2+y^2-2xy \cos \theta})(\sin\theta)^{2\alpha}\d\theta
\leq C.
$$
It follows that 
\begin{eqnarray}
\abs{T^\alpha_y\phi_j(x)}&\leq&C_N 2^{-2\alpha j}
\sum _{0\leq k<j} 2^{-k(N-2\alpha-1)}\chi_{\ent{(x-\frac{2^k}{x})_+,x+\frac{2^k}{x}}}(y)\nonumber\\
&&+C_N 2^{2(\alpha+1)j}\sum_{k\geq j}
2^{-kN}\chi_{\ent{(x-\frac{2^k}{x})_+,x+\frac{2^k}{x}}}(y),\label{eq:abc}
\end{eqnarray}
where $a_+=\max(0,a)$. Note that, since $2^{j-1}\leq x\leq 2^{j+1}$, $x-\frac{2^k}{x}\geq 0$ as long as
$k\leq 2j-2$.
From \eqref{eq:abc} we deduce that
\begin{eqnarray*}
\int_S  \abs{T^\alpha_y\phi_j(x)}\d \mu_\alpha(y)&\leq&
C_N 2^{-2\alpha j}\sum _{0\leq k<j} 2^{-k(N-2\alpha-1)}
\mu_\alpha\left(S\cap \ent{x-\frac{2^k}{x},x+\frac{2^k}{x}}\right)\\
&&+C_N 2^{2(\alpha+1)j}\sum_{j\leq k\leq 2j-2}
2^{-kN}\mu_\alpha\left(S\cap \ent{(x-\frac{2^k}{x})_+,x+\frac{2^k}{x}}\right)\\
&&+C_N 2^{2(\alpha+1)j}\sum_{k\geq 2j-1}
2^{-kN}\mu_\alpha\left(S\cap \ent{0,x+\frac{2^k}{x}}\right)\\
&=&C_N(\Sigma_1+\Sigma_2+\Sigma_3).
\end{eqnarray*}
Using \eqref{eq:trives1}, the first sum is simply estimated as follows:
$$
\Sigma_1\leq C2^{-2\alpha j}x^{2\alpha}\sum_{0\leq k<j} 2^{-k(N-2\alpha-2)} \eps
\leq C \sum _{k\ge 0} 2^{-k(N-2\alpha-2)}\eps\leq C\eps
$$
provided we take $N>2\alpha+2$.

For the second sum, we appeal again to \eqref{eq:trives1} and write
$$
\Sigma_2\leq C2^{(4\alpha+2)j}\sum_{j\le k<2 j} 2^{-k(N-1)}  \eps 
\leq C\eps
$$
provided we take $N>4\alpha+3$, while for the last sum we use \eqref{eq:trives2} to get
$$
\Sigma_3\leq C\frac{2^{2(\alpha+1)j}}{x^{2(\alpha+1)}} \sum_{ k\ge 2j} 2^{-kN}2^{2(\alpha+1)k}\eps
\leq C\eps.
$$
The proof of \eqref{eq:supy2} is similar.
\end{proof}
This completes the proof Theorem \ref{th:strong}.
\end{proof}
\begin{remark}
It would be interesting to obtain more precise quantitative estimates
of the constants $C(S,\Sigma)$
in Theorems \ref{weak} and \ref{th:strong}.
In a forthcoming work, we will obtain such an estimate in the case $S=[0,a]$ is an
interval and $\Sigma$ is $(\eps,\alpha)$-thin with $0<\eps<1$
arbitrary. This estimate takes the form $\norm{F_{\Sigma}E_{[0,a]}}\leq
f_a(\eps)$
where $f_a(\eps)\to0$ as $\eps\to0$.\footnote{An easy modification of
the above argument also provides such an estimate for $\eps$
\emph{small enough}.} Note that this allows to extend Theorem
\ref{th:strong}
to sets $S,\Sigma$ of the form $S=S_0\cup S_\infty$,
$\Sigma=\Sigma_0\cup\Sigma_\infty$ where
$S_0\subset [0,a]$, $\Sigma_0\subset[0,b]$ and
$S_\infty\subset[a,+\infty)$, $\Sigma_\infty\subset[b,+\infty)$
are $\eps$-thin.

Indeed, $F_\Sigma
E_S=F_{\Sigma_0}E_{S_0}+F_{\Sigma_\infty}E_{S_0}+F_{\Sigma_0}E_{S_\infty}+
F_{\Sigma_\infty}E_{S_\infty}$. Now, according to Theorem \ref{weak},
$\norm{F_{\Sigma_0}E_{S_0}}<1$.
Further, $\norm{F_{\Sigma_\infty}E_{S_0}}+\norm{F_{\Sigma_0}E_{S_\infty}}\leq
f_a(\eps)+f_b(\eps)\to0$
as $\eps\to 0$ and $\norm{F_{\Sigma_\infty}E_{S_\infty}}\leq
C\sqrt{\eps}$, according to (the proof of)
Theorem \ref{th:strong}. It follows that, if $\eps$ is small enough,
then $\norm{F_\Sigma E_S}<1$
so that $(S,\Sigma)$ is still a strong annihilating pair.
\end{remark}

\end{document}